\documentclass[a4paper,12pt]{article}

\usepackage[latin1]{inputenc}
\usepackage[T1]{fontenc}
\usepackage{latexsym,amssymb,amsmath,epsfig}

\newcommand{\sNN}{{\hbox{$\scriptstyle{I}$\kern-.25em\hbox{$\scriptstyle N$}}}}
\newcommand{\sRR}{{\hbox{$\scriptstyle{I}$\kern-.25em\hbox{$\scriptstyle R$}}}}
\newcommand{\sEE}{{\hbox{$\scriptstyle{I}$\kern-.25em\hbox{$\scriptstyle E$}}}}
\newcommand{\sZZ}{{\hbox{$\scriptstyle{Z}$\kern-.3em\hbox{$\scriptstyle Z$}}}}

\newtheorem{defi}{Definition}[section]

\newtheorem{prop}[defi]{Proposition}

\newtheorem{theo}[defi]{Theorem}
\newenvironment{dem}{\vskip 2mm\noindent {\it Proof}:}
		    {\hfill $\square$ \vskip 2mm \noindent}

\title{A zero-one law for first-order logic on random images}
\author{D. Coupier, A. Desolneux, and B. Ycart} 
\date{}
\begin{document}
\maketitle
\begin{center}
MAP5, FRE CNRS 2428, Universit\'e Ren\'e Descartes, Paris
\end{center}
\vskip 0.5cm\noindent
{\bf Mail address :} UFR Math-Info, 45 rue des Saints-P\`eres\\
\hspace*{3cm} 75270 PARIS CEDEX 06, FRANCE
\vskip 0.5cm\noindent
{\bf E-Mail address :} {\tt ycart@math-info.univ-paris5.fr}\\
{\bf Telephone :} 33 1 44 55 35 28\\
{\bf Fax :} 33 1 44 55 35 35
\vskip 2cm
\begin{abstract}
\noindent
For an $n\!\times\! n$ random image with independent pixels, black with
probability $p(n)$ and white with probability $1\!-\!p(n)$, 
the probability of satisfying any given
first-order sentence tends to $0$ or $1$, 
provided both $p(n)n^{\frac{2}{k}}$ and
$(1-p(n))n^{\frac{2}{k}}$ tend to $0$ or $+\infty$, for any
integer $k$. The result is proved by computing the threshold function 
for basic local sentences, and applying Gaifman's theorem.
\end{abstract}
\newpage
\setcounter{equation}{0}
\section{Introduction}
\label{intro}
The motivation for this work came for the Gestalt theory of vision
(see \cite{desolneux} and references therein), a basic idea of which
is that the human eye focuses first on remarkable or unusual
features of an image, i.e. features that would have a low probability
of occurring if the image were random. Hence the natural question:
which properties of a random image have a low or high probability?
Here we shall deal with the simplest model for random images:
\begin{defi}
\label{def:randimage}
Let $n$ be a positive integer. Consider the {\rm pixel set} 
$X_n=\{1,\ldots,n\}^2$.
An {\rm image of size $n\!\times\! n$} is a
mapping from $X_n$ to $\{0,1\}$
(white/black). Their set is denoted by $E_n$. It is endowed with the
product of $n^2$ independent copies of the Bernoulli distribution with
parameter $p$, that will be denoted by $\mu_{n,p}$~:
$$
\forall \eta\in E_n\;,\quad
\mu_{n,p}(\eta) = \prod_{i,j=1}^n p^{\eta(i,j)}(1-p)^{1-\eta(i,j)}\;. 
$$
A {\rm random image of size $n\!\times\! n$  and level $p$}, denoted by
${\cal I}_{n,p}$, is a random element of $E_n$ with distribution
$\mu_{n,p}$.  
\end{defi}
In other words, a random image of size $n\!\times\! n$ and level $p$ is a
square image in which all
pixels are independent, each being black with probability $p$ or white
with probability $1\!-\!p$.

We shall use the elementary definitions and concepts of first-order
logic on finite models, such as described for instance in Ebbinghaus 
and Flum \cite{EbbinghausFlum}. Gaifman's theorem 
(\cite{Gaifman82} and \cite{EbbinghausFlum} p.~31)
shows that first-order sentences are essentially local. They can be
logically reduced to the appearance of fixed subimages (precise
definitions will be given in section \ref{firstorder}). Assume $p$ is
fixed. Then as $n$ tends to infinity, any given subimage of fixed
size should appear somewhere in the random image ${\cal I}_{n,p}$,
with probability tending to $1$: this is
the two dimensional version of the well known ``typing monkey''
paradox. It justifies intuitively that the zero-one law should hold
for fixed values of $p$. Our main result is more general.
\begin{theo}
\label{th:principal}
Let $p(n)$ be a function from $\mathbb{N}$ into $[0,1]$ such that:
$$
\forall k=1,2,\ldots,\;
\lim_{n\rightarrow\infty} n^{\frac{2}{k}}p(n) = 0\mbox{ or } +\infty
\mbox{ and }
\lim_{n\rightarrow\infty} n^{\frac{2}{k}}(1-p(n)) = 0\mbox{ or } +\infty\;.
$$
Let $A$ be a first-order sentence. Then:
$$
\lim_{n\rightarrow\infty} 
Prob[\,{\cal I}_{n,p}\models A\,] = 0\mbox{ or } 1\;.
$$
\end{theo}
Zero-one laws have a long history (cf. Compton 
\cite{Compton} for a review and chapter 3 of \cite{EbbinghausFlum}). 
The first of them was proved
independently by Glebskii et al. \cite{glebskiietal} and Fagin
\cite{Fagin}. It applied to the first-order logic on a finite universe
without constraints, and uniform probability. As an example, 
interpret the elements of $E_n$ as directed graphs with vertex set 
$\{1,\ldots,n\}$, by putting an edge between $i$ and $j$ if pixel
$(i,j)$ is black. Then ${\cal I}_{n,p}$ becomes a random directed
graph (or digraph) with edge probability $p$ 
(see for instance \cite{Karp,Luczak}, or \cite{Bollobas} 
for a general reference). As a particular case of the Glebskii et al.
-- Fagin theorem, the zero-one law holds for first-order 
propositions on random digraphs. However, first-order logic 
on images is more expressive than on digraphs, since the 
geometry of images is not conserved in the graph interpretation. 

The theory of random (undirected) graphs was
inaugurated by Erd{\"o}s and R\'enyi
\cite{ErdosRenyi} (see \cite{Bollobas,Spencer} for general
references). The zero-one law holds for random graphs with edge
probability $p$, as a consequence of Oberschelp's theorem 
\cite{Oberschelp} on parametric classes (see \cite{EbbinghausFlum}
p.~74 or \cite{Spencer} p.~318). At first, zero-one laws were
essentially combinatorial, as they applied to the uniform probability
on the set of all structures, corresponding to edge probability
$p=\frac{1}{2}$ in the case of graphs. 
It was soon noticed that they also hold for any fixed 
value of $p$. But it is well known that random graphs become more 
interesting by letting $p=p(n)$ tend to $0$ as $n$ tends to
infinity. A crucial notion for random graphs
is the appearance of given subgraphs (\cite{Spencer}
p.~309). The threshold function for the
appearance of a given subgraph in a random graph is $p(n)=n^{-\frac{v}{e}}$,
where $v$ and $e$ are integers. For $p(n)=n^{-\frac{v}{e}}$, 
the probability of appearance for certain subgraphs does 
not tend to $0$ or $1$. Using the extension technique, 
(\cite{Gaifman64,Fagin} and
\cite{EbbinghausFlum} p.~73),
Shelah and Spencer \cite{ShelahSpencer} made a complete study of
those functions $p(n)$ for which the zero-one law holds for random
graphs, and proved in particular that it does for $p(n)=n^{-\alpha}$,
for any irrational $\alpha$. Theorem \ref{th:principal} is the analogue
for random images of Shelah and Spencer's result. To understand why,
first notice that the random image model is invariant through
exchanging black and white, together with $p$ and $1\!-\!p$. 
Thus we will consider only functions $p(n)$ tending to $0$. 
We shall define precisely the notion of threshold function in 
section \ref{threshold}, and prove that all threshold functions 
for patterns are of type $p(n)=n^{-\frac{2}{k}}$: the zero-one 
law does not hold for these values. For instance, if 
$p(n)$ is small (resp.: large) compared to $n^{-2}$, 
the probability of having at least one black pixel tends to
$0$ (resp.: $1$). But for $p(n)=n^{-2}$, it tends to $1-e^{-1}$. 
Theorem \ref{th:principal} essentially says that the zero-one 
law holds for any function $p(n)$ which is not a
threshold function.

It is worth pointing out here that theorem \ref{th:principal} can
be extended easily to other random structures, along two different
directions. Firstly, we chose to restrict the study to binary images,
using a single unary relation in the language (cf. section
\ref{firstorder}). With slight modifications of the proofs, and the
values of threshold functions, one could introduce a finite set of
 ``color'' unary relations, allowing for the coding of multilevel
gray or color images. The other possible generalisation concerns the
type of graphs. An image is essentially a colored square
lattice. The crucial property of that graph for our proof is that 
there exists a fixed number of vertices at fixed distance of 
any vertex (ball have bounded cardinality).  
Our study can easily be extended to any family of graphs with
bounded balls. For instance, theorem \ref{th:principal} also 
holds for a randomly colored $d$-dimensional square lattice
with $n^d$ points, up to replacing $n^{\frac{2}{k}}$ by 
$n^{-\frac{d}{k}}$ in its statement.

Section \ref{firstorder} is devoted to first-order logic on
images.  There we shall discuss basic local sentences
(definition \ref{def:basiclocal} and \cite{EbbinghausFlum} p.~31), and
reduce them to combinations of ``pattern sentences'' 
(definition \ref{def:pattern}), showing that a 
zero-one law holds for all first-order sentences if
it holds for basic local or pattern sentences 
(proposition \ref{prop:01patterns}). 
This will trivially imply that theorem \ref{th:principal} 
holds for fixed values of $p$. The section will
end with two examples of (second-order) sentences whose probability
under $\mu_{n,\frac{1}{2}}$ tends to $\frac{1}{2}$.

In section \ref{threshold}, we shall define the notion of threshold
function (definition \ref{def:threshold}) and prove that all 
threshold functions for basic local sentences are of type
$n^{-\frac{2}{k}}$ (proposition \ref{prop:basiclocal}). Theorem
\ref{th:principal} easily follows from propositions
\ref{prop:01patterns} and \ref{prop:basiclocal}.  
\setcounter{equation}{0}
\section{First-order logic for images}
\label{firstorder}
We shall follow the notations and definitions in chapter 0 of
\cite{EbbinghausFlum} for the syntax and semantics of first-order
logic. The {\it vocabulary} is the set of {\it relations} (or
predicates). They apply to the {\it universe} (or domain).
In our case the universe will be the pixel set $X_n$.
Image properties will not only be statements on colors of pixels
but also about their geometrical arrangement. Our vocabulary will
consist of $1$ unary and $4$ binary relations. The unary relation $C$ is
interpreted as the color: $Cx$ means that $x$ is a black pixel and
$\neg Cx$ that it is white. Before
defining the binary relations, we need a few considerations on the
geometry of $X_n$.

The pixel set $X_n$ is embedded in $\mathbb{Z}^2$, and naturally
endowed with a {\it graph structure}. In image analysis 
(see for instance chapter 6 of Serra \cite{Serra}), 
the cases most often considered are:
\begin{itemize}
\item
the {\it $4$-connectivity}. For $i,j>0$, the neighbors of $(i,j)$ are:
$$
(i+1,j), (i-1,j), (i,j+1), (i,j-1)\;.
$$   
\item
the {\it $8$-connectivity}. The $4$ diagonal neighbors are also
included:
$$
(i+1,j+1), (i-1,j+1), (i+1,j-1), (i-1,j-1)\;.
$$
\end{itemize}
At this point a few words about the borders are needed. In
order to avoid particular cases (pixels having less than $4$ or $8$ 
neighbors), we shall impose a periodic
boundary, deciding for instance that $(1,j)$ is neighbor with $(n,j), (n,j-1)$,
and $(n,j+1)$, so that the graph becomes a regular $2$-dimensional
torus. Although it may seem somewhat unnatural for images, without
that assumption the zero-one law would fail. Consider indeed
the (first-order) sentence ``there exist $4$ black pixels each having only one
horizontal neighbor''. Without periodic boundary conditions, 
it applies to the $4$ corners, and the probability for a random 
image ${\cal I}_{n,p}$ to satisfy it is $p^4$. From now on,
the identification $n+1\equiv 1$ holds for all operations on pixels. 

Once the graph structure is fixed, the relative positions of pixels
can be described by binary predicates. In the case of 
$4$-connectivity $2$ binary predicates suffice, $U$ (up) and $R$ (right):
$Uxy$ means that $y=x+(0,1)$ and $Rxy$ that $y=x+(1,0)$. In the case of 
$8$-connectivity, two more predicates must be added, $D_1$ and $D_2$:
$D_1xy$ means that $y=x+(1,1)$ and $D_2xy$ that $y=x+(1,-\!1)$. 
For convenience reasons, we shall stick to $8$-connectivity. Thus the
vocabulary of images is the set $\{C,U,R,D_1,D_2\}$. Once the
universe and the vocabulary are fixed, the {\it structures} are particular
models of the relations, applied to variables in the domain. To any
structure, a graph is naturally associated (\cite{EbbinghausFlum} p.~26),
connecting those pairs of elements $\{x,y\}$ which are such that
$Sxy$ or $Syx$ are satisfied, where $S$ is any of the binary relations.
Of course only those structures for which the associated graph 
is the square lattice with diagonals and periodic boundaries 
will be called images. As usual, the {\it graph distance}
$d$ is defined as the minimal length of a path between two
pixels. We shall denote by $B(x,r)$ the ball of center $x$ and radius
$r$:
$$
B(x,r) = \{\, y\in X_n\,;\; d(x,y)\leq r\,\}\;
$$
In the case of $8$-connectivity, $B(x,r)$ is a square containing 
$(2r+1)^2$ pixels.
\vskip 3mm
Formulas such as $Cx$, $Uxy$, $Rxy$\ldots are called {\it atoms}. The 
{\it first-order logic} (\cite{EbbinghausFlum} p.~5)
is the set of all formulas obtained by recursively combining 
first-order formulas, starting with atoms.
\begin{defi}
The set ${\cal L}_1$ of first-order formulas is defined by:
\begin{enumerate}
\item
All atoms belong to ${\cal L}_1$.
\item
If $A$ and $B$ are first-order formulas, then
$(\neg A)$, $(\forall x Ax)$ and $(A\wedge B)$ also belong to 
${\cal L}_1$.
\end{enumerate}
\end{defi}
Here are two examples of first-order formulas:
\begin{enumerate}
\item
$\forall x,y,\;(Rxy\wedge Uyz)\rightarrow D_1xz$,
\item
$(\exists y\;(Rxy\wedge Uyz)) \leftrightarrow D_1xz$
\end{enumerate}
Notice that any image satisfies them both: adding the
two diagonal relations $D_1$ and $D_2$ does not make the language any 
more expressive. The only reason why the $8$-connectivity was preferred
here is that the corresponding balls are squares.

\begin{figure}[!ht]
\centerline{
\includegraphics[width=6cm,height=6cm]{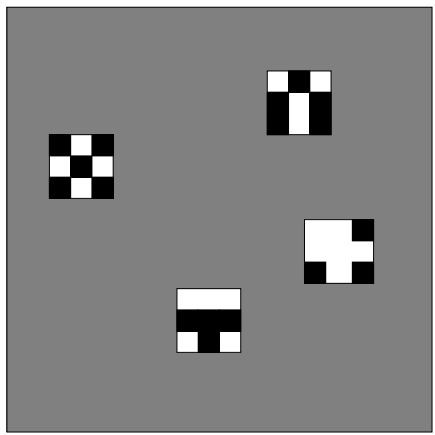}
} 
\caption{Illustration of a pattern sentence, for $m=4$ and $r=1$.}
\label{fig:pattern}
\end{figure}

\begin{figure}[!ht]
\centerline{
\includegraphics[width=6cm,height=6cm]{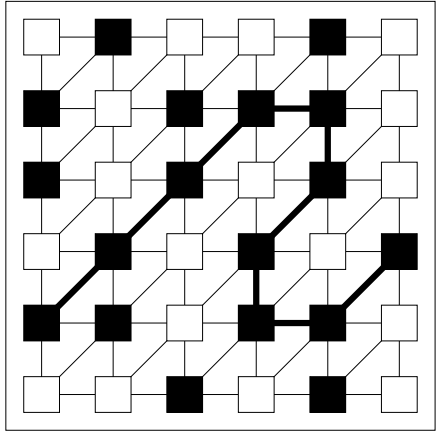}
} 
\caption{A $6$-connected path of black pixels from left to right.}
\label{fig:connectedpath}
\end{figure}

We are interested in formulas for which it can be decided if they are
true or false for any given image, i.e. for which all variables are
quantified. They are called {\it closed formulas}, or {\it
sentences}. Such a sentence $A$ defines a subset $A_n$ of $E_n$: that of
all images $\eta$ that {\it satisfy} $A$ ($\eta\models A$). Its
probability for $\mu_{n,p}$ will still be denoted by $\mu_{n,p}(A)$.
$$
\mu_{n,p}(A) = Prob[ {\cal I}_{n,p} \models A ]
= \sum_{\eta\models A}\mu_{n,p}(\eta)\;.
$$
Gaifman's theorem (\cite{EbbinghausFlum} p.~31), states that every
first-order sentence is equivalent to a boolean combination of 
{\it basic local sentences}.
\begin{defi}
\label{def:basiclocal}
A {\it basic local sentence} has the form:
\begin{equation}
\label{deflocalebasique}
\exists x_1\ldots\exists x_m\; 
\left(\bigwedge_{1\leq i<j\leq m} d(x_i,x_j)>2r\right)
\;\wedge\;
\left(\bigwedge_{1\leq i \leq m} \psi_i(x_i) \right)\;,
\end{equation}  
where: 
\begin{itemize}
\item[$\bullet$] $m$ and $r$ are fixed nonnegative integers,
\item[$\bullet$] for all $i=1,\ldots,m$, 
$\psi_i(x) \in {\cal L}_1$ is a formula for which only variable $x$ is
free (not bound by a quantifier), and the other variables
all belong to the ball $B(x,r)$.
\end{itemize}
\end{defi}
For any $x$ and a fixed radius $r$, 
consider now a {\it complete description} 
$D(x)$ of the ball $B(x,r)$, i.e. a first-order sentence for which 
all statements concerning pixels at distance at most $r$ of $x$ are 
either asserted or negated. There exists a single image $I_D$ of size 
$(2r+1)\!\times\! (2r+1)$, centered at $x$, satisfying it. Thus $D(x,r)$ 
can be interpreted as: ``the pattern of
pixels at distance at most $r$ of $x$ is $I_D$''. 
\begin{defi}
A {\it pattern sentence} has the form:
\label{def:pattern}
\begin{equation}
\label{defpattern}
\exists x_1\ldots\exists x_m\; 
\left(\bigwedge_{1\leq i<j\leq m} d(x_i,x_j)>2r\right)
\;\wedge\;
\left(\bigwedge_{1\leq i \leq m} D_i(x_i) \right)\;,
\end{equation}  
where: 
\begin{itemize}
\item[$\bullet$] $m$ and $r$ are fixed nonnegative integers,
\item[$\bullet$] for all $i=1,\ldots,m$, 
$D_i(x)$ is a complete description of the ball $B(x,r)$.
\end{itemize}
\end{defi}
Examples of (interpreted) pattern sentences are:
\begin{enumerate}
\item
``there exist $3$ black pixels'',
\item
``there exists a $3\!\times\! 3$ white square'',
\item
``there exist $3$ non overlapping $5\!\times\! 5$ white squares 
with a black pixel on the center''.
\end{enumerate}
Figure \ref{fig:pattern} gives another illustration.
Obviously, pattern sentences are particular cases of basic
local sentences. Proposition \ref{prop:01patterns} below reduces the
proof of zero-one laws for random images to pattern sentences. 
\begin{prop}
\label{prop:01patterns}
Consider the following three assertions.
\begin{itemize}
\item[(i)] The probability of any {\rm pattern} sentence tends to $0$ or $1$.
\item[(ii)] The probability of any {\rm basic local} sentence 
tends to $0$ or $1$.
\item[(iii)] The probability of any {\rm first order} sentence 
tends to $0$ or $1$.
\end{itemize}
Then {\it (i)} implies {\it (ii)} and {\it (ii)} implies {\it (iii)}.
\end{prop}
\begin{dem}
Observe first that if the probabilities of sentences $A$ and $B$ tend
to $0$ or $1$, then so do the probabilities of $\neg A$ and $A\wedge
B$. This follows from elementary properties of probabilities.
As a consequence, if the probability of $A$ tends to $0$ or $1$ 
for any $A$ in a
given family, this remains true for any finite boolean combination of
sentences in that family. Thus Gaifman's theorem yields that
{\it (ii)} implies {\it (iii)}. 
We shall prove now that every basic local sentence is either
unsatisfiable or a finite 
boolean combination of pattern sentences. Indeed, consider a formula 
$\psi(x)$ for which only variable $x$ is free, and the other variables
all belong to the ball $B(x,r)$. Either it is not satisfiable, 
or there exists a finite set of
$(2r+1)\!\times\!(2r+1)$ images (at most $2^{(2r+1)^2}$) which satisfy
it. To each of these images corresponds a complete description
$D(x)$ which implies $\psi(x)$. So $\psi(x)$ is equivalent 
to the disjunction of these $D(x)$'s:
\begin{equation}
\label{decompositionlocalcomplete}
\psi(x) \leftrightarrow \bigvee_{D(x)\rightarrow \psi(x)} D(x)\;.
\end{equation}
In formula (\ref{deflocalebasique}), one can replace each $\psi_i(x_i)$ 
by a disjunction of complete descriptions. Rearranging terms, one sees
that the basic local sentence (\ref{deflocalebasique}) is itself a
finite disjunction of pattern sentences.
\end{dem}
The zero-one law for fixed values of $p$ is an easy consequence of
proposition \ref{prop:01patterns}. Indeed, for fixed $p$, the
probability of any pattern sentence tends to $1$. 
To see why, consider the following
sentence:
\begin{equation}
\label{inlinepattern}
\exists x\;
\left(\bigwedge_{1\leq i \leq m} D_i(x+((i-1)(2r+1),0)) \right)\;,
\end{equation}
interpreted as: ``subimages $I_{D_1},\ldots,I_{D_m}$ appear in $m$
consecutive, horizontally adjacent balls of radius $r$''. 
It clearly implies (\ref{defpattern}). But (\ref{inlinepattern}) is
equivalent to the appearance of a given subimage on a rectangle of
size $(2r+1)\!\times\! (m(2r+1))$. This occurs in a random image 
${\cal I}_{n,p}$ with probability tending to $1$ as $n$ tends to
infinity. Thus (\ref{defpattern}) has a probability tending to one of being
satisfied by ${\cal I}_{n,p}$.
\vskip 3mm
This section ends with two counter-examples of (second-order)
sentences the probability of which does not tend to $0$ or $1$. The
first one is ``the number of black pixels is even''. Its probability
is $\frac{1}{2}(1+(1-2p)^{n^2})$, which tends to $\frac{1}{2}$ for any $p$
such that $0<p<1$. The second example is more visual. Define a 
{\it $6$-connected path} as an $m$-tuple of pixels $(x_1,\ldots,x_m)$,
such that for $i=1,\ldots,m\!-\!1$, $x_{i+1}\in
x_i\pm\{(1,0),(0,1),(1,1)\}$, and the borders of the image are not 
crossed (see an illustration on figure \ref{fig:connectedpath}). 
Consider now the two sentences:
\begin{enumerate}
\item $BLR$: ``there exists a $6$-connected path of black pixels from
left to right'',
\item $WTB$: ``there exists a $6$-connected path of white pixels from
top to bottom''. 
\end{enumerate}
Some geometrical considerations show that an image satisfies $BLR$ if
and only if it does not satisfy $WTB$ (this would not hold for $4$- or
$8$-connected paths: see \cite{Serra} p.~183). Take now
$p=\frac{1}{2}$. Symmetry implies that 
$\mu_{n,\frac{1}{2}}(BLR) = \mu_{n,\frac{1}{2}}(WTB)$. 
Hence both probabilities must be equal to $\frac{1}{2}$.

The sentences $BLR$ and $WTB$ are examples of those properties 
studied by percolation theory (see Grimmett 
\cite{Grimmett} for a general reference). Actually the random image
model that we consider here is a finite approximation of site
percolation (\cite{Grimmett} p.~24). 
Using percolation techniques, one can prove that $\mu_{n,p}(BLR)$
tends to $0$ if $p<\frac{1}{2}$, to $1$ if $p>\frac{1}{2}$. 
\setcounter{equation}{0}
\section{Threshold functions for basic local sentences}
\label{threshold}
The notions studied in this section have exact counterparts in the
theory of random graphs as presented by Spencer
\cite{Spencer}.  
We begin with the asymptotic probability of single
pattern sentences, which correspond to the appearance of subgraphs
(\cite{Spencer} p.~309).
\begin{prop}
\label{prop:thresholdpattern}
Let $r$ and $k$ be two integers such that $0<k<(2r+1)^2$.
Let $I$ be a fixed $(2r+1)\!\times\! (2r+1)$ image, with $k$ black pixels
and $h=(2r+1)^2-k$ white pixels. Let $D(x)$ be the complete
description of the ball $B(x,r)$ satisfied only by a copy of $I$,
centered at $x$. Let $\tilde{D}$ be the sentence $(\exists x\; D(x))$.
Let $p=p(n)$ be a function from $\mathbb{N}$ to $[0,1]$.
\begin{equation}
\label{threshold1}
\mbox{If } \lim_{n\rightarrow\infty} n^2p(n)^k =0
\mbox{ then } \lim_{n\rightarrow\infty} \mu_{n,p(n)}(\tilde{D}) = 0\;.
\end{equation}
\begin{equation}
\label{threshold2}
\mbox{If } \lim_{n\rightarrow\infty} n^2p(n)^k(1-p(n))^h =+\infty
\mbox{ then } \lim_{n\rightarrow\infty} \mu_{n,p(n)}(\tilde{D}) = 1\;.
\end{equation}
\begin{equation}
\label{threshold3}
\mbox{If } \lim_{n\rightarrow\infty} n^2(1-p(n))^h =0
\mbox{ then } \lim_{n\rightarrow\infty} \mu_{n,p(n)}(\tilde{D}) = 0\;.
\end{equation}
\end{prop}
\begin{dem}
We already noticed the symmetry of the problem: switching black and
white together with $p$ and $1\!-\!p$ should leave statements
unchanged. In particular the proofs of (\ref{threshold1}) and
(\ref{threshold3}) are symmetric, and only the former will be given.

For a given $x$, the probability of occurence of $I$ in the ball
$B(x,r)$ is~:
$$
\mu_{n,p(n)}(D(x))  = p(n)^{k}(1-p(n))^h\;.
$$
The pattern sentence $\tilde{D}$ is the disjunction of all $D(x)$'s:
$$
\tilde{D} \leftrightarrow \bigvee_{x\in X_n} D(x)\;.
$$
Hence:
$$
\mu_{n,p(n)}(\tilde{D})\leq n^2 p(n)^{k}(1-p(n))^h\;,
$$
from which (\ref{threshold1}) follows. 

Consider now the following set
of pixels:
\begin{equation}
\label{deftiling}
T_n = \{\,(r+1+\alpha(2r+1),r+1+\beta(2r+1))\,,\;\alpha,\beta=0,\ldots,
{\scriptstyle \lfloor\frac{n}{2r+1}\rfloor}\!-\!1\,\}\;, 
\end{equation}
where $\lfloor\,\cdot\,\rfloor$ denotes the integer part. 
Call $\tau(n)$ the cardinality of $T(n)$:
$$
\tau(n) = \left\lfloor\frac{n}{2r+1}\right\rfloor^2\;,
$$
which is of order $n^2$. Notice that the disjunction of $D(x)'s$ for
$x\in T_n$ implies $\tilde{D}$.
$$
 \bigvee_{x\in T_n} D(x) \rightarrow \tilde{D}\;.
$$
The distance between any two distinct pixels $x,y \in T_n$ is larger 
than $2r$, and the balls $B(x,r)$ and $B(y,r)$ do not overlap. 
Therefore the events 
``${\cal I}_{n,p}\models D(x)$'' for $x\in T_n$ are mutually independent.
Thus:
\begin{eqnarray*}
\mu_{n,p(n)}(\tilde{D})&\geq&
\mu_{n,p(n)}\left(\bigvee_{x\in T_n}D(x)\right)\\
&=& 1 - \Big(1-p(n)^{k}(1-p(n))^h\Big)^{\tau(n)}\\
&\geq& 1-\exp(-\tau(n)p(n)^k(1-p(n))^h)\;,
\end{eqnarray*}
hence (\ref{threshold2}).
\end{dem}
Due to the symmetry of the model, we shall consider from now on that
$p(n)<\frac{1}{2}$. Proposition \ref{prop:thresholdpattern} shows that
the appearance of a given subimage only depends on its number of black
pixels: if $p(n)$ is small compared to $n^{-\frac{2}{k}}$, then no
subimage of fixed size, with $k$ black pixels, should appear in 
${\cal I}(n,p(n))$. If $p(n)$ is large compared to $n^{-\frac{2}{k}}$,
all subimages with $k$ black pixels should appear.
Proposition \ref{prop:thresholdpattern} does not cover the particular
cases $k=0$ (appearance of a white square) and $k=(2r+1)^2$ (black
square). They are easy to deal with. Denote by $W$ (resp.: $B$) the
pattern sentence $(\exists x\; D(x))$, where $D(x)$ denotes the complete
description of $B(x,r)$ being all white (resp.: all black). Then 
$\mu_{n,p(n)}(W)$ always tends to $1$ (remember that
$p(n)<\frac{1}{2}$). Statements (\ref{threshold1}) and 
(\ref{threshold2}) apply to $B$, with $k=(2r+1)^2$. 
\vskip 3mm
The notion of {\it threshold function} is a formalisation of the behaviors
that have just been described.
\begin{defi}
\label{def:threshold}
Let $A$ be a sentence. A threshold function for $A$ is a
function $r(n)$ such that:
$$
\lim_{n\rightarrow\infty} \frac{p(n)}{r(n)} = 0
\mbox{ implies }
\lim_{n\rightarrow\infty}
\mu_{n,p(n)}(A) = 0\;,
$$
and~:
$$
\lim_{n\rightarrow\infty} \frac{p(n)}{r(n)} = +\infty
\mbox{ implies }
\lim_{n\rightarrow\infty}
\mu_{n,p(n)}(A) = 1\;.
$$
\end{defi}
Notice that a threshold function is not unique. For instance if
$r(n)$ is a threshold function for $A$, then so is $cr(n)$ for any
positive constant $c$. It is costumary to ignore this and talk about
``the'' threshold function of $A$. For instance, the threshold
function for ``there exists a black pixel'' is $n^{-2}$.

Proposition \ref{prop:thresholdpattern} essentially says that the
threshold function for the appearance of a given subimage $I$ is
$n^{-\frac{2}{k}}$, where $k$ is the number of black pixels in $I$.
Proposition \ref{prop:basiclocal} below will show that the threshold
function for a basic local sentence $L$ is $n^{-\frac{2}{k(L)}}$, where
$k(L)$ is an integer that we call the {\it index} of $L$. 
Its definition uses the decomposition
(\ref{decompositionlocalcomplete}) of a local property 
into a finite disjunction of complete
descriptions, already used
in the proof of proposition \ref{prop:01patterns}.
\begin{defi}
\label{def:index}
Let $L$ be the basic local sentence defined by:
$$
\exists x_1\ldots\exists x_m\; 
\left(\bigwedge_{1\leq i<j\leq m} d(x_i,x_j)>2r\right)
\;\wedge\;
\left(\bigwedge_{1\leq i \leq m} \psi_i(x_i) \right)\;.
$$
If $L$ is not satisfiable, then we shall set $k(L)=+\infty$.
If $L$ is satisfiable,
for each $i=1,\ldots,m$, consider the finite set 
$\{ D_{i,1},\ldots,D_{i,d_i}\}$ of
those complete descriptions on the ball $B(x_i,r)$ which imply
$\psi_i(x_i)$. 
$$
\psi_i(x_i) \leftrightarrow \bigvee_{1\leq j\leq d_i}
D_{i,j}(x_i)\;.
$$
Each complete description $D_{i,j}(x_i)$ corresponds to an image on
$B(x_i,r)$. Denote by $k_{i,j}$ its number of black pixels.

The {\rm index} of $L$, denoted by $k(L)$ is defined by:
\begin{equation}
\label{defindex}
k(L) = \max_{i=1}^m \min_{j=1}^{d_i} k_{i,j}\;. 
\end{equation}
\end{defi}
The intuition behind definition \ref{def:index} is the following. Assume 
$p(n)$ is small compared to $n^{-\frac{2}{k(L)}}$.  Then there exists 
$i$ such that none of the $D_{i,j}(x_i)$ can be satisfied, therefore 
there is no $x_i$ such that $\psi_i(x_i)$ is satisfied, and $L$ is not
satisfied. On the
contrary, if $p(n)$ is large compared to $n^{-\frac{2}{k(L)}}$, then
for all $i=1,\ldots,m$, $\psi_i(x_i)$ should be satisfied for at least
one pixel $x_i$, and the
probability of satisfying $L$ should be large. In other words,
$n^{-\frac{2}{k(L)}}$ is the threshold function of $L$. 
\begin{prop}
\label{prop:basiclocal}
Let $L$ be a basic local property, and $k(L)$ be its index. If $L$ is
satisfiable and $k(L)>0$, then its threshold function is 
$n^{-\frac{2}{k(L)}}$. If $k(L)=0$, its probability tends to $1$
(for $p(n)<\frac{1}{2}$).
\end{prop}
\begin{dem}
Assume $L$ is satisfiable (otherwise its probability is null) and $k(L)>0$.
Let $r(n)=n^{-\frac{2}{k(L)}}$.
For $p(n)<\frac{1}{2}$, we need to prove that $\mu_{n,p(n)}$ tends to
$0$ if $p(n)/r(n)$ tends to $0$, and that it tends to $1$ if
$p(n)/r(n)$ tends to $+\infty$. The former will be proved first.

Consider again the decomposition of $L$ into complete descriptions:
$$
L \leftrightarrow 
\exists x_1\ldots\exists x_m\; 
\left(\bigwedge_{1\leq i<j\leq m} d(x_i,x_j)>2r\right)
\;\wedge\;
\left(\bigwedge_{1\leq i \leq m} \bigvee_{1\leq j\leq d_i}
D_{i,j}(x_i)\right)\;.
$$
If $p(n)/r(n)$ tends to $0$, there exists $i$ such that:
$$
\forall j=1,\ldots, d_i\,,\;
\lim_{n\rightarrow\infty} n^2p(n)^{k_{i,j}} = 0\;.
$$
By proposition \ref{prop:thresholdpattern}, the probability of 
$(\exists x\;D_{i,j}(x))$ tends to zero for all $j=1,\ldots,d_i$. 
Therefore the probability of 
$(\exists x\,\psi_i(x))$ tends to $0$, which implies that
$\mu_{n,p(n)}(L)$ tends to $0$.

Conversely, for each $i=1,\ldots,m$, choose one of the $D_{i,j}(x)$'s, 
such that the number of black pixels in the corresponding image is
minimal (among all $k_{i,j}$'s). Denote that particular description by
$D_i(x)$. Consider now the following pattern sentence, which implies
$L$:
\begin{equation}
\label{pattern1}
\exists x_1\ldots\exists x_m\; 
\left(\bigwedge_{1\leq i<j\leq m} d(x_i,x_j)>2r\right)
\;\wedge\;
\left(\bigwedge_{1\leq i \leq m}
D_i(x_i)\right)\;.
\end{equation}
As in the proof of proposition \ref{prop:thresholdpattern}, 
we shall use the lattice $T_n$, defined by
(\ref{deftiling}). Remember that its cardinality $\tau(n)$ is of order
$n^2$. The pattern sentence (\ref{pattern1}) is implied by:
\begin{equation}
\label{pattern2}
\exists x_1\ldots\exists x_m\; 
\left(\bigwedge_{1\leq i\leq m} x_i\in T_n\right)
\;\wedge\;
\left(\bigwedge_{1\leq i<j\leq m} x_i\neq x_j \right)
\;\wedge\;
\left(\bigwedge_{1\leq i \leq m}
D_i(x_i)\right)\;.
\end{equation}
Assume first that $k(L)=0$. Then necessarily, for each $i$, the image
corresponding to $D_i(x)$ has only white pixels. With
$p(n)<\frac{1}{2}$, the probability of observing a
$(2r+1)\!\times\!(2r+1)$ white image is larger than $\pi=2^{-(2r+1)^2}$. Since
subimages centered at the points of $T_n$ are independent, the
probability of (\ref{pattern2}) is larger than:
$$
1-\sum_{l=0}^{m-1} \binom{\tau(n)}{l}\pi^{l}(1-\pi)^{\tau(n)-l}\;,
$$
which tends to $1$ as $n$ tends to infinity.

Assume now that $k(L)>0$. The images corresponding to the minimal
descriptions $D_i$ need not be all different: 
renumber {\it different} descriptions $D_i$ as
$D'_1,\ldots, D'_{m'}$. Let $k(i)$ be the number
of black pixels of $D'_i$ (hence $k(L)=\max\{k(i)\}$).
Let $\pi_i(n)$ be the probability of $D'_i(x)$, for a given $x$:
$$
\pi_i(n) = p(n)^{k(i)}(1-p(n))^{(2r+1)^2-k(i)}\;.
$$
>From the random image ${\cal I}_{n,p}$ define
the random variable $N_i$ as the number of those pixels $x_i\in T_n$ such
that ${\cal I}_{n,p}$ is described by $D'_i(x_i)$ on the ball
$B(x_i,r)$. Since the different balls do not overlap, 
$N_i$ has a binomial distribution, with parameters $\tau(n)$ and
$\pi_i(n)$. Assuming $p(n)/r(n)$ tends to $+\infty$, it is easy to
check that the product $\tau(n)\pi_i(n)$ also tends to infinity.
Therefore the probability that $N_i$ is larger than $m$ tends to $1$,
and so does the probability that {\it all} the $N_i$'s are larger than
$m$. But if all the $N_i$'s are larger than $m$, then
${\cal I}_{n,p}$ satisfies (\ref{pattern2}), hence (\ref{pattern1}) and
$L$.
\end{dem}
Having characterized the threshold functions of all basic local
properties, the proof of theorem \ref{th:principal} is now clear.
If $p(n)n^{\frac{2}{k}}$ tends to $0$ or $+\infty$ for any positive 
integer $k$, then by proposition
\ref{prop:basiclocal} the probability of any basic local sentence 
tends to $0$ or $1$. This remains true for any boolean combination of
basic local sentences (cf. proposition \ref{prop:01patterns}). By
Gaifman's theorem, these boolean combinations cover all first-order  
sentences. Hence the zero-one law for first-order logic.
\bibliography{zero_un}
\bibliographystyle{plain}
\end{document}